# Optimization of the Matching Network for using Genetic Algorithm

Jalil Rasekhi[1], Jalil Rashed Mohasel[2]

[1]j.rasekhi@ece.ut.ac.ir

[2]jrashed@ut.ac.ir

*Abstract* — **Microstrip-like antenna (MLA) which was developed nearly a decade ago, is a powerful radiating element. The primary challenge in designing a MLA is to provide an optimized matching network such that the overall input reflection is kept as low as possible within the required bandwidth. In this paper, the necessity and procedure of applying genetic algorithm to MLA problems has been presented. Comparison with the existing literature shows good agreement in the overall input reflection.**

## I. INTRODUCTION

Development trends of antennas with conformal topology, light weight, superior performance, low cost, easy maintenance and modular design have been on the move for decades. This has led to many kinds of antennas currently in use. One of these superior performance antennas is microstrip-like antenna or simply MLA. This antenna which was originated and developed in TUdelft [1], is a very economical microwave solution, since it is just a dielectric filled waveguide. It is not only low-cost and low-profile but at the same time enjoys a very small cross-polarization compared to alternative microstrip antennas. On the other hand, the MLA antennas face a crucial obstacle which is their high aperture reflection. Hajian and Tran applied air-gap matching technique (AMT) to match MLA and integrate it with antenna circuitry[2]. The technique requires finding optimal values for matching network. In that work a try and error method was used. Here we present a genetic solution to carry out the problem solving. In the next sections the necessity and application of GA method for this case will be introduced.

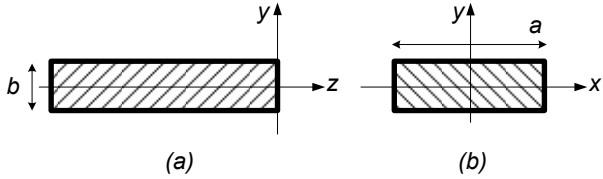

Fig. 1. MLA antenna: (a) Side view and (b) Front view

## II. PROBLEM DESCRIPTION

A schematic of an MLA antenna is depicted in Fig.1.
The most important characteristic of this antenna is that its largest dimension is less than half of the free space wavelength. In addition to this advantage, MLA could provide another benefit if used in an array. Since they are smaller than half a wavelength, MLA's may be arranged in array grids so that no grating-lobe effects would appear[3]. However, as mentioned in section 1, MLA has such a high aperture reflection which prevents any practical use. Tran [2] presented a technique to match this high input reflection using an air-gap and two E-plane steps. MLA and its matching network are illustrated in Fig.2. In this figure, the values of $l_i$ and $b_i$ must be found such that the overall input reflection at reference plane $P_1$ be minimized.

Fig. 2. MLA and its matching network (side view)

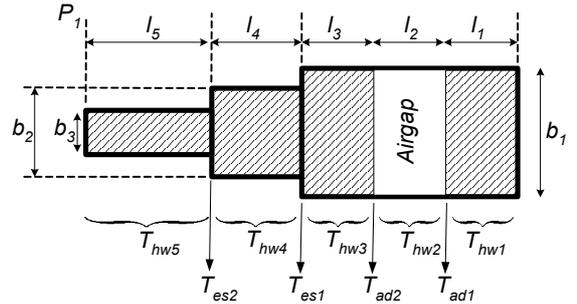

To start with, let the transmission matrices of different sections of the matching structure be defined as:

$T_{es}$ : E-plane step transmission matrix
$T_{hw}$ : Homogeneous waveguide transmission matrix
$T_{ad}$ : Airgap-to-dielectric transmission matrix

Formulas for these matrices are readily available in [2]. Therefore the transmission matrix of the overall matching network can be found as:

$$T = \begin{bmatrix} T_{11} & T_{12} \\ T_{21} & T_{22} \end{bmatrix}$$

$$= T_{hw1} T_{ad1} T_{hw2} T_{ad2} T_{hw3} T_{es1} T_{hw4} T_{es2} T_{hw5} \quad (1)$$

Hence the input reflection coefficient at reference plane $P_1$, may be written as [2]:

$$G_{in} = \frac{T_{11} G_{ap} + T_{12}}{T_{21} G_{ap} + T_{22}} \quad (2)$$

where, $G_{ap} = \dfrac{1 - y_{ap}}{1 + y_{ap}}$, is MLA aperture reflection coefficient. The evaluation procedure for the aperture admittance $y_{ap}$ is also available in [2]:

$$y_{ap} = \frac{Y_{11}}{Y_{10}} + 2\sum_{m=3,5,\ldots} D_m \frac{Y_{m1}}{Y_{10}}$$

$$+ \sum_{m=3,5,\ldots} D_m^2 \cdot \left(\frac{Y_{mm}}{Y_{10}} + \frac{Y_{m0}}{Y_{10}}\right) \quad (3)$$

where

$$Y_{ij} = Y_{ji} = \frac{2}{ab}\frac{1}{wm}\frac{1}{4\pi^2} \times$$

$$\int_{-\infty}^{+\infty}\int_{-\infty}^{+\infty} \frac{(k^2 - k_x^2)}{k_z} C_0^2(k_y) C_i(k_x) C_j(k_x) dk_y dk_x$$

$$D_m = -\frac{Y_{m1}}{Y_{mm} + Y_{m0}}$$

$$C_0(k_y) = \frac{b\sin(\frac{k_y b}{2})}{\frac{k_y b}{2}}$$

$$C_m(k_x) = \frac{2m\pi a j^{m-1}\cos(\frac{k_x a}{2})}{(m\pi)^2 - (k_x a)^2}$$

This integral is evaluated using Parseval's theorem and calculus of residues. Detailed discussions regarding the integration procedure are available in [4] and [5],[10]. The aperture reflection is then calculated by considering three modes (modes 1, 3 and 5).

As mentioned earlier, our goal is to find optimal values for $l_i$ and $b_i$, such that the resultant $G_{in}$ is minimized. The most recent work used a try and error method with good results[2]. We make use of a genetic approach to carry out the procedure.

## III. MOTIVES FOR EMPLOYING GENETIC OPTIMIZATION FOR THE MLA PROBLEM

It should be noticed that MLA problem has an essential difference with other optimization problems which is the unspecified number of air-gaps and E-plane steps from the beginning. Hence we can not find a definite general cost function. Even if a pre-defined matching structure is assumed, cost function could be very complex and should be re-evaluated for any changes in the structure[2],[8],[9]. On the other hand, if we use genetic algorithm we just need to evaluate transmission matrices of different MLA sections and carry out the multiplication in equation (1). In this way we not only avoid difficulties and complexity of finding an equation for the overall input reflection, but the programming procedure will also be more simple and may be easily developed to include more general cases. This can not be achieved by using other optimization methods.

## IV. OPTIMIZATION ALGORITHM[6], [7]

The utilized optimization flowchart is illustrated in Fig.3.

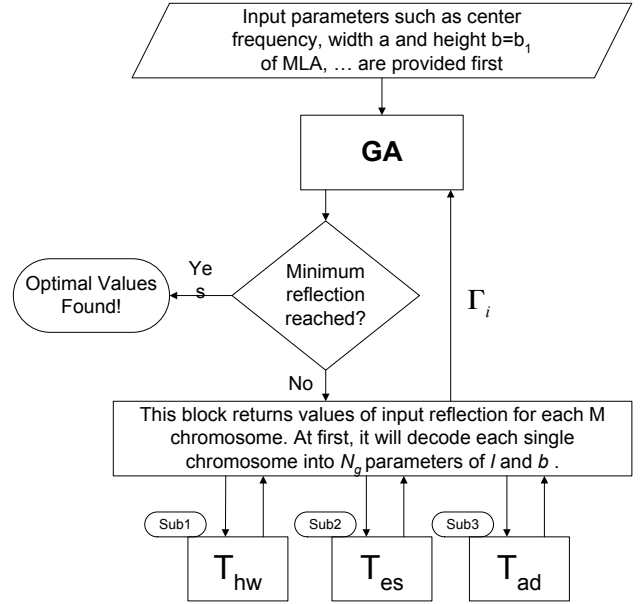

Fig. 3. Optimization flowchart

GA is the main body of this flowchart and includes genetic processes such as mutation, reproduction, natural selection and so on. When routine is initialized, GA will produce $M$ random chromosomes containing $N_g$ 8-bit genes. These genes are base 2 encoded versions of optimization values $l_i$ and $b_i$. Each chromosome then, has $N = 8N_g$ bits (see figure 4). After initial chromosome generation, GA calls for $M$ values of $G_i$ relating to these $M$ chromosomes.

"Mlaref.m" subroutine is responsible for evaluating $G_i$ s and performs this in three steps:

- Decodes $M$ received chromosomes into $M$ configuration vectors:
  $$Mlacfg^i = \{l_1^i, l_2^i, \ldots, l_5^i, b_1^i, b_2^i, b_3^i\}$$
- Referring to subroutines "$Sub1$", "$Sub2$" and "$Sub3$", evaluates transmission matrices of each $Mlacfg^i$ and multiplies them to obtain overall transmission matrix of the matching network for $i$ th chromosome.
- Determines $G_{in}^i$ for these $M$ transmission matrices using equation (2)

The procedure is repeated for younger generations until a global minimum is reached by $\min(G_{in}^i)$.

Fig. 4. Structure of chromosomes and the way genes are organized in each chromosome

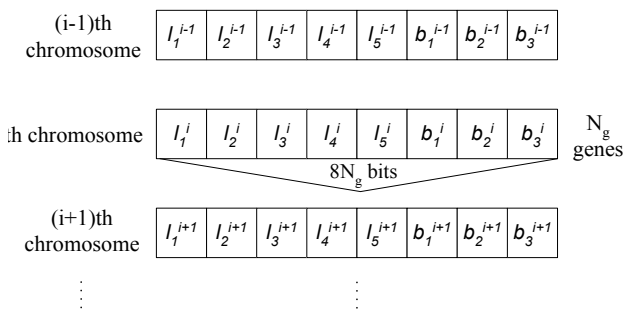

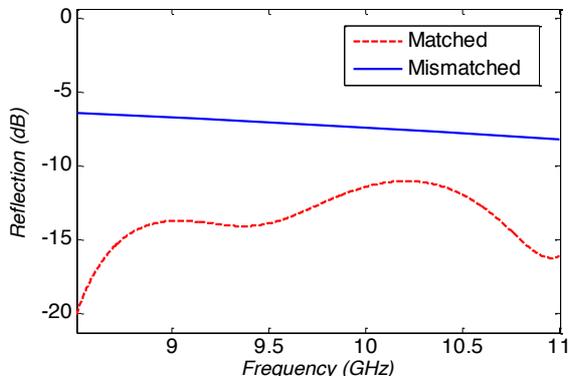

Fig. 5. Resulting input reflection compared to that of mismatched case Results

## V. RESULTS

The above illustrated genetic algorithm was applied to an MLA antenna with dimensions of $a = 17mm$, $b = 11mm$ with a centre frequency of $9.75 GHz$. A set of optimized parameters $\{l_1, l_2, ... l_5, b_1, b_2, b_3\}$ were found. Simulation results with these matching network parameters showed good results and improvement for input reflection compared to that of [2]. Fig. 5 illustrates the resultant input reflection employing genetically optimized parameters. The simulations were carried out using HPHFSS package.

As seen from the same figure, wideband matching has been possible for MLA and input reflection is more assuring on the two extremes of the design band.

## REFERENCES

[1] M. Tian, P.D. Tran, M. Hajian and L.P. Ligthart, Air-gap technique for matching the aperture of miniature waveguide antennas, IEEE instrumentation and measurement technology conference, IMtc/93, Irvine, California, May 1993.

[2] M. Hajian, D.P. Tran and L.P. Ligthart, Design of a wideband miniature dielectric-filled waveguide antenna for collision-avoidance radar, IEEE Trans. AP, Vol. 42, No. 1, February 2000.

[3] D.M. Pozar and D.H. Schaubert, Analysis of an infinite array of rectangular microstrip patches with idealized probe feeds, IEEE trans. AP, Vol. AP-32, No. 10, Oct. 1984.

[4] M. Tian, "Characterization of miniature Dielectric Filled Open Ended Waveguide Antennas", PhD thesis, Delft University, October 1995.

[5] J. Rasekhi, "Design and analysis of a wideband Miniaturized waveguide antenna in X-band," M.Sc. thesis, Tehran University; September 2004.

[6] R. L. Haupt, An Introduction to Genetic Algorithms for Electromagnetics, IEEE Antenna and Propagation Magazine, Vol. 37, No. 2, April 1995.

[7] J. H. Holland, "Genetic Algorithms", Scientific American, July 1992, pp. 66-72.

[8] Rastegar, S.; Babaeian, A.; Bandarabadi, M.; Toopchi, Y., "Airplane detection and tracking using wavelet features and SVM classifier," in *System Theory, 2009. SSST 2009. 41st Southeastern Symposium on*, vol., no., pp.64-67, 15-17 March 2009

[9] Babaeean, A.; Tashk, A.B.; Barzin, F.; Hosseini, S.M., "Target Tracking Using Mean Shift and Dynamic Directional Gradient Vector Flow," in *System Theory, 2008. SSST 2008. 40th Southeastern Symposium on*, vol., no., pp.366-370, 16-18 March 2008

[10] Tashk, A.R.B.; Faez, K., "Boosted Bayesian Kernel Classifier Method for Face Detection," in *Natural Computation, 2007. ICNC 2007. Third International Conference on*, vol.1, no., pp.533-537, 24-27 Aug. 2007